\documentclass[11pt]{article}

\usepackage{amsfonts}
\usepackage{epsfig,graphics,subfigure,psfrag,amsmath,amssymb}
\usepackage{latexsym,wrapfig,picinpar,picins}
\usepackage{bm}
\usepackage{amsthm}
\usepackage{graphicx}
\usepackage{mathrsfs}
\usepackage{fancyref}
\usepackage{dcolumn}
\usepackage{color}
\usepackage{pstricks}
\usepackage{epstopdf}
\usepackage{subfigure}
\usepackage{multicol}
\usepackage{enumerate}
\usepackage{caption}
\usepackage{geometry}
\usepackage{threeparttable}
\makeatletter
\@addtoreset{equation}{section}
\makeatother

\newtheorem{thm}{Theorem}[section]

\newtheorem{lem}[thm]{Lemma}

\def\E\mathbb{ E}

\newcommand{\be}{\begin{eqnarray}}
\newcommand{\ee}{\end{eqnarray}}
\newcommand{\ben}{\begin{eqnarray*}}
\newcommand{\een}{\end{eqnarray*}}

\begin{document}

\title{Convergence analysis of a numerical scheme for the porous medium equation by an energetic variational approach}

\author{Chenghua Duan\footnotemark[2]\and Chun Liu\footnotemark[3]\and Cheng Wang\footnotemark[4]\and Xingye Yue\footnotemark[1]}

 \renewcommand{\thefootnote}{\fnsymbol{footnote}}
 \footnotetext[2]{Department of Mathematics, Soochow University, Suzhou 215006, China and Shanghai Center for Mathematical Sciences,  Fudan University, Shanghai 200438, China (chduan@fudan.edu.cn).}
 \footnotetext[3]{Department of Applied Mathematics, Illinois Institute of Technology, Chicago, IL 60616, USA,
 (cliu124@iit.edu).}

 \footnotetext[4]{Department of Mathematics, University of Massachusetts, Dartmouth, North Dartmouth, MA, 02747-2300, USA (cwang1@umassd.edu).}
 \footnotetext[1]{Corresponding author. Department of Mathematics, Soochow University, Suzhou 215006, China (xyyue@suda.edu.cn).}
%



\date{}
\maketitle

{
\begin{abstract}
\footnotesize
The porous medium equation (PME) is a typical nonlinear degenerate parabolic equation. We have studied numerical methods for PME by an energetic variational approach in [C. Duan et al, J. Comput. Phys., 385 (2019) 13-32], where  the trajectory equation can be obtained and two numerical schemes have been developed based on different dissipative energy laws.   It is also proved that the nonlinear scheme, based on  $f \log f$ as the total energy form of the dissipative law, is uniquely solvable on an admissible convex set and preserves the corresponding discrete dissipation law. Moreover, under certain smoothness assumption, we have also obtained the second order convergence in space and the first order convergence in time for the scheme. In this paper, we provide a rigorous proof of the error estimate by a careful higher order asymptotic expansion and two step error estimates. The latter technique contains a rough estimate to control the highly nonlinear term in a discrete $W^{1,\infty}$ norm, and a refined estimate is applied to derive the optimal error order.

 \end{abstract}
}

\section{Introduction and Background }
\label{sec:1}
 One of the typical nonlinear degenerate parabolic equations is the porous medium equation (PME):  $$ \partial_t f=\Delta_{x} (f^{m}), \ x\in\Omega\subset\mathbb{R}^d,\  m>1,$$
where  $f:=f(x,t)$ is a non-negative scalar function of space $x\in\mathbb{R}^d$ ($d \ge 1$) and the time $t\in\mathbb{R}^{+}$, and $m$ is a constant larger than 1.  It has been applied in many physical and biological models, such as an isentropic gas flow through a porous medium,  the viscous gravity currents, nonlinear heat transfer and image processing \cite{J. L. Vazquez(2007)}, etc.  

It is well known that the PME is  degenerate at points where $f = 0$.
In turn, the PME has many special features: the finite speed of  propagation, the free boundary, a possible waiting time phenomenon \cite{C. Duan(2018),J. L. Vazquez(2007)}.
Various numerical methods have been studied for the PME, such as  finite difference approach \cite{J. L. Graveleau(1971)}, tracking algorithm method \cite{E. DiBenedetto(1984)},  a  local discontinuous Galerkin finite element method \cite{Q. Zhang(2009)},  Variational Particle Scheme (VPS) \cite{M. Westdickenberg(2010)} and  an adaptive moving mesh finite element method \cite{C. Ngo(2017)}. Many theoretical analyses have been derived in the existing literature \cite{D.G. Aronson(1969), A.S. Kalasnikov(1967), O.A. Oleinik(1958), S. Shmarev(2003),S. Shmarev(2005), J. L. Vazquez(2007)}, etc.

Relevant detailed descriptions can be found in a recent paper \cite{C. Duan(2018)}, in which the numerical methods for the PME were constructed by an Energetic Variational Approach (EnVarA) to naturally keep the physical laws, such as the conservation of mass, energy dissipation and force balance. Meanwhile, based on different dissipative energy laws, two different numerical schemes have been studied. In more details, based on  the total energy form  $f \log f$ and  $\frac{1}{2 f}$, a fully discrete nonlinear scheme and a linear numerical  scheme could be appropriately designed for the trajectory equation, respectively. It has also been proved that  the former one is uniquely solvable on an admissible  convex set, and both schemes preserve the corresponding discrete dissipation law. Numerical experiments  have demonstrated that both schemes have yielded a good approximation for the solution without oscillation and the free boundary.  The notable advantage is that the waiting time problem could be naturally treated, which has been a well-known difficult issue for all the existing methods. In addition, under certain smoothness assumption, the second order convergence in space and the first order convergence in time have been reported for both schemes in \cite{C. Duan(2018)}. The aim of the paper is to provide a rigorous proof of the optimal rate convergence analysis for the nonlinear scheme. On the other hand, the highly nonlinear nature of  the trajectory equation makes the convergence analysis every challenging. To overcome these difficulties, we use a higher order expansion technique to ensure a higher order consistency estimate, which is needed to obtain a discrete $W^{1,\infty}$ bound of the numerical solution. Similar ideas have been reported in earlier literature for incompressible fluid equations~\cite{W. E(1995), W. E(2002), R. Samelson(2003), C. Wang(2000)}, while the analysis presented in this work turns out to be more complicated, due to the lack of a linear diffusion term in the trajectory equation of the PME. In addition, we have to carry out two step estimates to recover the nonlinear analysis:
\begin{itemize}
\item {\bf Step 1} A rough estimate for the discrete derivative of numerical solution, namely ($D_h x^{n+1}_h$)  at time $t_{n+1}$, to control the nonlinear term;
\item {\bf Step 2} A refined estimate for the numerical error function to obtain an optimal convergence order.
\end{itemize}
Different from a standard error estimate, the rough estimate controls the nonlinear term, which is an effective approach to handle the highly nonlinear term.

   This paper is organized as follows. The trajectory equation of the PME and  the numerical scheme are outlined in Sec. \ref{sec:2_1} and  Sec. \ref{sec:2_2}, respectively. Subsequently, the proof of  optimal rate convergence analysis is provided in Sec. \ref{sec:3}. Finally we present a simple numerical example to demonstrate the convergence rate of the numerical scheme in Sec. \ref{sec:4}.

\section{Trajectory equation and the numerical scheme}
In this section, we  review the trajectory equation and the corresponding numerical scheme. 

\subsection{Trajectory equation of the PME}
\label{sec:2_1}
  In this part,  the one-dimensional trajectory equation will be reviewed, derived by an Energetic Variational Approach~\cite{C. Duan(2018)}. We solve  the following initial-boundary problem:
\begin{eqnarray}
&&
  \partial_t f+\partial_{x} (f{\bf v})=0, \ x\in\Omega\subset\mathbb{R}, \ t>0,
 \label{eqcm}\\
&& f{\bf v}=-\partial_{x}(f^{m}), \ x\in\Omega, \  m>1, \label{eqDa}\\
&& f(x,0)=f_{0}(x)\geq 0, \ x\in\Omega, \label{eqini} \\
&&\partial_{x}f=0, \ x\in\partial\Omega,\ t>0, 
 \label{eqbc}
 \end{eqnarray}
where $\Omega$ is a bounded domain,   $f$ is a non-negative function,   $t$ is the time, $x$ is the particle position and  ${\bf v}$ is the velocity of particle.

The following lemma is available.

\begin{lem} \label{lem-dissip}
  $f(x,t)$ is a positive solution of \eqref{eqcm}-\eqref{eqbc} if and only if
  $f(x,t)$ satisfies the corresponding energy dissipation law:
  \begin{equation}\label{equ:energylaw}
   \frac{d}{dt}\int_{\Omega}f\ln fdx=-\int_{\Omega}\frac{f}{mf^{m-1}}|{\bf v}|^{2}dx.
  \end{equation}
    \end{lem}

\noindent\textbf{Proof}:
 We first prove the energy dissipation law (\ref{equ:energylaw}) if $f$ is the solution of  \eqref{eqcm}-\eqref{eqbc}. Multiplying  by $(1+\ln f)$  and integrating on both sides of (\ref{eqcm}), we get
$$\int_{\Omega}(1+\ln f)\partial_t fdx=-\int_{\Omega}(1+\ln f)\partial_{x} ( f{\bf v})dx.$$
Using integration by parts, in combination with \eqref{eqDa}, we have
\begin{equation}\label{3}
  \frac{d}{dt}\int_{\Omega}f\ln fdx=\int_{\Omega}\frac{\partial_x f}{f} (f{\bf v})dx\\
=-\int_{\Omega}\frac{f}{mf^{m-1}}|\textbf{v}|^{2}dx\leq 0.
\end{equation}

Subsequently, we are also able to derive  \eqref{eqDa}  from the energy dissipation law  (\ref{equ:energylaw}) by EnVarA.

 In addition, (\ref{eqcm}) is the conservation law. In the Lagrangian coordinate, its solution can be expressed by:
 \begin{equation}\label{equ:conservationL}
  f(x(X,t),t)=\frac{f_{0}(X)}{\frac{\partial x(X,t)}{\partial X}},
\end{equation}
where $f_{0}(X)$ is the positive initial data and $\partial_X x:=\frac{\partial x(X,t)}{\partial X}$ is the \emph{deformation\ gradient} in one dimension.
Based on an Energetic Variational Approach, we can obtain the trajectory equation.
\begin{itemize}
\item Energy Dissipation Law.

  The {\bf total energy} of the PME is
  \begin{equation} \label{tot-energ-e} E^{total}:=\int_{\Omega}f\ln fdx. \end{equation}
\item Least Action Principle step. \ \

With (\ref{equ:conservationL}), the action functional in Lagrangian coordinate becomes
$$\mathcal{A}(x):=\int^{T^{*}}_{0}(-\mathcal{H}) dt = -\int^{T^{*}}_{0}
\int_{\Omega}f_0(X)\ln \Big(\frac{f_0(X)}{\partial_X x}\Big)dXdt,$$
where  $T^{*}>0$ is a given  terminal time and  $\mathcal{H}$ is the free energy depending on $x$.
Thus for any test function $y(X,t)=\widetilde{y}(x(X,t),t)\in C_{0}^{\infty}(\Omega\times (0,T^*))$ and $\varepsilon\in \mathbb{R}$, taking the variational of $\mathcal{A}(x)$ with respect to $x$, we have
\begin{subequations}
  \begin{align}
  \frac{d}{d\varepsilon}\bigg|_{\varepsilon=0}\mathcal{A}(x+\varepsilon y)%
&=\int^{T^{*}}_{0}\int_{\Omega}\frac{f_0(X)}{\partial_X x}\cdot\partial_X y\ dXdt \nonumber\\
  &=-\int^{T^{*}}_{0}\int_{\Omega}\partial_{x}f\cdot\widetilde{y}\  dxdt \nonumber.
  \end{align}
\end{subequations}
  Then the conservation force turns out to be
  $$F_{con}= \frac{\delta\mathcal{A}}{\delta x} =-\partial_{x}f,$$ in the Eulerian coordinate, and
   $$F_{con}=-\partial_{X}\left(\frac{f_0(X)}{\partial_X x}\right),$$ in the Lagrangian coordinate.

\item  Maximal Dissipation Principle step.  \ \

  Define the entropy production
$\Delta:=\int_{\Omega}\frac{f}{mf^{m-1}}|{\bf v}|^{2}dx$. Taking the variational of $\frac 12 \Delta$ with respect to the velocity ${\bf v}$  and ${x}_{t}$, we obtain the dissipation force
$$F_{dis}:=\frac{\delta\frac{1}{2}\Delta}{\delta{\bf v}}=\frac{f}{mf^{m-1}}{\bf v},$$ in the Eulerian coordinate
and
$$F_{dis}:=\frac{\delta\frac{1}{2}\Delta}{\delta (\partial_{t} {x})}=\frac{f_{0}(X)}{m\big(\frac{f_0(X)}{\partial_X x}\big)^{m-1}}\partial_{t} {x},$$ in the Lagrangian coordinate.\\
\item Force  balance step. \ \ Based on the Newton's force balance law, we get \end{itemize}

\begin{equation} \label{trajLa}
\frac{f_{0}(X)}{m\big(\frac{f_0(X)}{\partial_X x}\big)^{m-1}}\partial_{t} {x}=-\partial_{X}\left(\frac{f_0(X)}{\partial_X x}\right),
\end{equation}
in the Lagrangian coordinate, and the Darcy's Law in the Eulerian coordinate
\begin{equation*}
    \begin{aligned}
    \frac{f}{mf^{m-1}}{\bf v}= -\partial_{x} f.
    \end{aligned}
\end{equation*}
$\hfill\Box$\\


It is noticed that, there is an assumption that the value of initial state $f_0(x)$ is positive  in $\Omega$ to  make $\int_{\Omega}f\ln fdx$ well-defined in \eqref{equ:energylaw}. More details can be found in \cite{C. Duan(2018)}.


Since then,  we should first settle the initial and boundary conditions for \eqref{trajLa}. From   \eqref{eqcm} and \eqref{eqbc}, we have $ x_t|_{\partial\Omega}=0$, for $t>0$. This means that the particles lying on boundary will stay there forever, so  a Dirichlet boundary condition should be subject to as  $x|_{\partial\Omega}= X|_{\partial\Omega}$, for $t\ge 0$. As a result, the trajectory problem becomes


 \begin{eqnarray}
 &&    \frac{f_{0}(X)}{m\big(\frac{f_0(X)}{\partial_X x}\big)^{m-1}}\partial_{t} {x}=-\partial_{X}\left(\frac{f_0(X)}{\partial_X x}\right),\ X\in\Omega,\ t>0,  \label{eqtra}\\
 && x|_{\partial\Omega}= X|_{\partial\Omega},\ t>0,\label{eqtrabou} \\
     && x(X,0)=X, \ X\in\Omega. \label{eqtraini}
    \end{eqnarray}
Finally, with a substitution of \eqref{eqtra} into \eqref{equ:conservationL}, we obtain the solution $f(x,t)$ to \eqref{eqcm}-\eqref{eqbc}.

\subsection{Numerical scheme of the trajectory equation}
\label{sec:2_2}
 Let $\tau=\frac{T}{N}$, where $N\in\mathbb{N}^+$,  $T$  is the final time, and the grid points are given by $t_n=n\tau$, $n=0,\cdots,N$. Let $X_0$ be the left point of $\Omega$ and $h=\frac{|\Omega|}{M}$  be the spatial step,  $M\in\mathbb{N}^{+}$. Denote by $X_{r}=X(r)=X_0+ r h$, where $r$ takes on integer and half integer values.  Let $\mathcal{E}_{M}$ and $\mathcal{C}_{M}$ be the spaces of functions whose domains are $\{X_{i}\ |\ i=0,...,M\}$ and $\{X_{i-\frac{1}{2}}\ |\ i=1,...,M\}$, respectively. In component form, these functions are identified via
$l_{i}=l(X_{i})$, $i=0,...,M$, for $l\in\mathcal{E}_{M}$, and $\phi_{i-\frac{1}{2}}=\phi(X_{i-\frac{1}{2}})$,  $i=1,...,M$, for $\phi\in\mathcal{C}_{M}$.

\indent{The} difference operator $D_{h}: \mathcal{E}_{M}\rightarrow\mathcal{C}_{M}$, $d_{h}: \mathcal{C}_{M}\rightarrow\mathcal{E}_{M}$, and $\widetilde{D}_h: \mathcal{E}_{M}\rightarrow\mathcal{E}_{M}$  can be defined as:
\begin{align}\label{equ:dif1}
& (D_{h}l)_{i-\frac{1}{2}}= (l_{i}-l_{i-1})/h,\ i=1,...,M, \\
&  (d_{h}\phi)_{i}= (\phi_{i+\frac{1}{2}}-\phi_{i-\frac{1}{2}})/h,\ i=1,...,M-1, \\
& (\widetilde{D}_hl)_{i}= (l_{i+1}-l_{i-1})/2h,\ i=1,...,M-1, \label{equ:dif2_1}\\
& (\widetilde{D}_hl)_{i}= (4l_{i+1}-l_{i+2}-3l_{i})/2h,\ i=0, \label{equ:dif2_2} \\
& (\widetilde{D}_hl)_{i}= (l_{i-2}-4l_{i-1}+3l_{i})/2h,\ i=M,\label{equ:dif2_3}
\end{align}
respectively.

Let $\mathcal{Q}:=\{l \in\mathcal{E}_{M}\ |\ l_{i-1}<l_{i},\ 1\leq i\leq M;\ l_{0}=X_0,\ l_{M}=X_M\}$ with its boundary set $\partial\mathcal{Q}:=\{l \in\mathcal{E}_{M}\ |\ l_{i-1}\leq l_i,\ 1\leq i\leq M;\ \ l_{0}=X_0,\ l_{M}=X_M;   \exists\, i\in \{1,\ldots, M\},\ s.t.\  l_{i-1}=l_{i}\}$. Then $\bar{\mathcal{Q}}:=\mathcal{Q}\cup\partial\mathcal{Q}$ is a closed convex set. Its physical meaning indicates that particles are arranged in the order without twisting or exchanging in $\mathcal{Q}$.

A few more notations have to be introduced.  Let $l$, $g\in\mathcal{E}_M$ and $\phi$, $\varphi\in\mathcal{C}_M$. We define the \emph{inner product} on space $\mathcal{E}_M$ and $\mathcal{C}_M$ respectively as:
\begin{eqnarray}
  &&
 \left\langle l , g \right\rangle  := h \left( \frac12 l_0 g_0 + \sum_{i=1}^{M-1}   l_{i} g_{i}
  + \frac12 l_M g_M \right),\\
  \label{FD-inner product-1}
&&
  \left\langle \phi , \varphi \right\rangle_e
  := h \sum_{i=0}^{M-1}   \phi_{i+\frac 12} \varphi_{i+\frac 12} .
  \label{FD-inner product-2}
\end{eqnarray}

The following summation by parts formula is available:
\begin{equation}
   \left\langle l , d_h \phi \right\rangle
  = -\left\langle D_h l , \phi \right\rangle_e,
  \mbox{\ with $l_0 = l_M =0$, $\phi\in\mathcal{C}_{M}$, $l\in\mathcal{E}_M$}.
  \label{FD-inner product-3-1}
\end{equation}

The inverse inequality is given by:
\begin{equation}\label{equ:inverse}
\|l\|_{\infty}\leq C_m \frac{\|l\|_2}{h^{1/2}}, \ \ \forall l\in\mathcal{E}_M, \quad
\mbox{with $\|l\|_{\infty}:=\max\limits_{0\leq i\leq M} \{l_i\}\mbox{\ \ and\ \ }\|l\|_2^2:=\left\langle l,l \right\rangle$} .
\end{equation}

 \indent The {\bf  fully discrete scheme} is formulated as follows: Given the positive initial state $f_0(X)\in\mathcal{E}_M$ and the particle position $x^{n} \in\mathcal{Q}$, find $x^{n+1}=(x^{n+1}_{0},...,x^{n+1}_{M})\in\mathcal{Q}$ such that
 \begin{equation}\label{equ:numnum0}
\frac{f_{0}(X_{i})}{m\big(\frac{f_0(X)}{\widetilde{D}_h x^{n}}\big)_{i}^{m-1}}\cdot\frac{x^{n+1}_{i}-x^{n}_{i}}{\tau}=-d_{h}\Big[\Big(\frac{f_{0}(X)}{D_{h} x^{n+1}}\Big)\Big]_{i}, \ 1\leq i \leq M-1,
\end{equation}
 with $x_0^{n+1}=X_0$ and $x_M^{n+1}=X_M$, $n=0,\cdots,N-1$.

It is noticed that \eqref{equ:numnum0} is still a nonlinear system which can be solved by Newton's iteration method  \cite{C. Duan(2018)}. Then we obtain the numerical solution $f(x_i,t^n):=f^n_i$  by
\begin{equation}\label{eqDen}
f^n_i=\frac{f_0(X)}{\widetilde{D}_h x^n_i},\ 0\leq i\leq M,
\end{equation}
which is the discrete scheme of \eqref{equ:conservationL}.

\section{Convergence analysis}
\label{sec:3}
In this section,  the second order spatial convergence and the first order temporal convergence will be theoretically justified for the numerical scheme \eqref{equ:numnum0}.
We first introduce a higher order  approximate expansion  of the exact solution, since a consistency estimate (second order in space and first order in time) is not able to control the discrete $W^{1,\infty}$ norm of the numerical solution. Also see the related works in the earlier literature~\cite{A. Baskaran(2013), W. E(1995), W. E(2002), Z. Guan(2017), Z. Guan(2014), R. Samelson(2003), C. Wang(2004), C. Wang(2002), C. Wang(2000), L. Wang(2015)}, etc.

\begin{lem}
Assume a higher order  approximate solution  of the exact solution $x_e$:
\begin{equation}\label{approximation}
W:=x_e+\tau w_\tau^{(1)}+\tau^2 w^{(2)}_\tau +h^2 w_h,
\end{equation}
where $w_\tau^{(1)}$, $w_\tau^{(2)}$, $w_h \in C^{\infty}{(\Omega; 0,T)}$. Then there exists a small $\tau_0>0$, such that $\forall\tau, h \leq \tau_0$, $\widetilde{D}_h W>0$, i.e., $W\in\mathcal{Q}$, where $\tau$ and  $h$ are the time step and the spatial mesh sizes, respectively.
\end{lem}
\noindent{\textbf{Proof:}}
Because of a point-wise condition for the exact solution, $x_e\in\mathcal{Q}$, i.e.,  $\exists$ $\varepsilon_0>0$, such that $D_h x_e>\varepsilon_0>0$.
For small $\tau_0$, such that $\|\tau D_h w_\tau^{(1)}\|_{L^{\infty}}\leq\frac{1}{9}\varepsilon_0$, $\|\tau^2 D_h w_\tau^{(2)}\|_{L^{\infty}}\leq\frac{1}{9}\varepsilon_0$ and $\|h^2 D_h w_h\|_{L^{\infty}}\leq\frac{1}{9}\varepsilon_0$, for $\forall\tau, h \leq \tau_0$.
As a consequence, for $\forall \tau, h\leq\tau_0$, we have
\begin{equation}\label{U_Q}
D_h W\geq\frac{1}{3}\varepsilon_0>0 ,
\end{equation}
which in turn implies that $W\in\mathcal{Q}$.   $\hfill\Box$

\begin{thm}
\label{convergence}
Assume that  the initial function $f_0(X)$ is positive and bounded, i.e., $0<b_f \leq f_0(X) \leq B_f$.
 Denote $x_e\in\Omega$ as the exact solution to the original PDE \eqref{eqtra} (with enough regularity) and $x_{h}\in\mathcal{Q}$ as the numerical solution to \eqref{equ:numnum0}. The numerical error function is defined at a point-wise level:
\be\label{error function-1}
 e_i^n = x_{e_i}^n - x_{h_i}^n ,
\ee
where $x_{e_i}^n,\ x_{h_i}^n\in\mathcal{Q}$, $0 \le i \le N$, $n=0,\cdots,M$. 

Then we have
\begin{itemize}
\item{
$e^{n}=(e^{n}_0,\cdots,e^{n}_M)$ satisfies}
$$\quad  \| e^{n}\|_2:=\langle e^{n},e^{n}\rangle\le C (\tau + h^2).$$
\item{$\widetilde{D}_h e^{n}=(\widetilde{D}_h e^{n}_0,\cdots,\widetilde{D}_h e^{n}_M)$ satisfies}
$$\quad  \| \widetilde{D}_h e^{n} \|_{2} \leq C(\tau+h^2).$$
Moreover, the error between the numerical solution $f_h^{n}$ and the exact solution $f_e^{n}$ of equaiton \eqref{eqcm}-\eqref{eqbc}  can be estimated by:
$$\quad  \|f_h^{n}-f_e^{n}\|_2 \leq C(\tau+h^2),$$
where $C$ is a positive constant, $h$ is the spatial step, $\tau$ is the time step and $n=0,\cdots,N$.
\end{itemize}
\end{thm}

\noindent\textbf{Proof}:
A careful Taylor expansion of the exact solution in both time and space, in terms of the numerical scheme~\eqref{eqtra}, gives that
\begin{eqnarray}
&&
 \frac{f_{0}(X_{i})}{m(\frac{f_0(X_{i})}{\widetilde{D}_h x_{e_{i}}^{n}})^{m-1}}\frac{x_{e_{i}}^{n+1}-x_{e_{i}}^{n}}{\tau}=-d_{h}\Big(\frac{f_{0}(X)}{D_{h} x_{e}^{n+1}}\Big)_{i}+ \tau l^{(1)}_i +\tau^{2}l^{(2)}_i+\tau^{3}l^{(3)}_i+h^{2}g^{(1)}_i+h^4g^{(2)}_i,\nonumber\\
&&
\ \ \ \ \ 1\leq i\leq M-1,\nonumber\\
&&
 \quad \mbox{with $x_{e_{0}}^{n+1} = X_0$ ,   $\quad  x_{e_{M}}^{n+1} = X_M $}, \label{consistency-2}
\end{eqnarray}
where $\|l^{(1)}\|_2$, $\|l^{(2)}\|_2$, $\|l^{(3)}\|_2$, $\|g^{(1)}\|_2$, $\|g^{(2)}\|_2\leq C_e$,  with $C_e$ only dependent on the exact solution.

To perform a higher order consistency analysis for an approximate solution  of the exact solution, we have to construct the approximation $W$ as in \eqref{approximation}.

The term $w_\tau^{(1)}\in C^{\infty}{(\Omega; 0,T)}$ is given by the following linear equation:
\begin{eqnarray}
  && \frac{f_0(X)}{m(\frac{f_0(X)}{\partial_X x_e})^{m-1}}\partial_t w_\tau^{(1)}+  \frac{m-1}{m (\frac{f_0(X)}{\partial_X x_e})^{m-2}} \partial_t x_e\cdot\partial_X w_\tau^{(1)}
  =\partial_X\Big(\frac{f_0(X)}{(\partial_X x_e)^2}\partial_X w_\tau^{(1)}\Big)-l^{(1)}, \label{equ:u_1}\notag\\
  &&
 w_\tau^{(1)}|_{\partial\Omega}=0,  \label{u_1_bou} \ \ \ \   w_\tau^{(1)}(\cdot,0)=0.
  \label{u_1_ini}
\end{eqnarray}

The term $w^{(2)}_\tau\in C^{\infty}{(\Omega; 0,T)}$ is given by the following linear equation:
\begin{eqnarray}
  && \frac{f_0(X)}{m(\frac{f_0(X)}{\partial_X x_e})^{m-1}}\partial_t w^{(2)}_\tau+  \frac{m-1}{m (\frac{f_0(X)}{\partial_X x_e})^{m-2}} \partial_t x_e\cdot\partial_X w^{(2)}_\tau \nonumber\\
  &&
  +\frac{(m-1)(m-2)}{2m (\frac{f_0(X)}{\partial_X x_e})^{m-2}\partial_X x_e}(\partial_X w^{(1)}_\tau)^2\cdot\partial_t x_e +\frac{(m-1)}{m (\frac{f_0(X)}{\partial_X x_e})^{m-2}}\partial_t w^{(1)}_\tau\cdot \partial_X w^{(1)}_\tau\nonumber \\
  &&
=\partial_X\Big(\frac{f_0(X)}{(\partial_X x_e)^2}\partial_X w^{(2)}_\tau\Big)-\partial_X\Big(\frac{f_0(X)}{(\partial_X x_e)^3}(\partial_X w^{(1)}_\tau)^2\Big)-l^{(2)},\label{equ:u_2}\notag \\
  &&
  w^{(2)}_\tau|_{\partial\Omega}=0,  \label{u_2_bou} \ \ \ \  w^{(2)}_\tau(\cdot,0)=0.
  \label{u_2_ini}
\end{eqnarray}

The term $w_h\in C^{\infty}{(\Omega; 0,T)}$ is given by the following linear equation:
\begin{eqnarray}
  &&
  \frac{f_0(X)}{m(\frac{f_0(X)}{\partial_X x_e})^{m-1}}\partial_t w_h+  \frac{(m-1)\partial_t x_e}{m (\frac{f_0(X)}{\partial_X x_e})^{m-2}} \partial_X w_h=\partial_X\Big(\frac{f_0(X)}{(\partial_X x_e)^2}\partial_X w_h\Big)-g^{(1)}, \label{equ:u_h}\notag \\
  &&
  w_h|_{\partial\Omega}=0, \label{u_h_bou}\ \ \ \   w_h(\cdot,0)=0.
  \label{u_h_ini}
\end{eqnarray}
Since  $w_\tau^{(1)}$, $w_\tau^{(2)}$, $w_h$ are dependent only on $W$ and $x_e$, we have the following estimate:
\begin{eqnarray}
  &&
  \|W-x_e\|_{H^m}=\tau\|w_\tau^{(1)}\|_{H^m}+\tau^2\|w_\tau^{(1)}\|_{H^m}+h^2\|w_h\|_{H^m}\leq C'(\tau+h^2). \label{estimate_u1}
  %
\end{eqnarray}

With these expansion terms, the constructed approximation $W\in\mathcal{Q}$ satisfies the numerical scheme with a higher order truncation error:
\begin{eqnarray}\label{num:U}
&&
\frac{f_{0}(X_i)}{m\big(\frac{f_0(X)}{\widetilde{D}_h W^{n}}\big)_i^{m-1}}\cdot\frac{W^{n+1}_i-W^{n}_i}{\tau}=-d_{h}\Big(\frac{f_{0}(X)}{D_{h} W^{n+1}}\Big)_i+ \tau^{3}l^{*}_i+h^4g^{*}_i, \ \ 1\leq i\leq M-1, \nonumber\\
&&
\mbox{with\ \ } W_0^{n+1}=X_0, \ \ W_M^{n+1}=X_M, \ \ n=0,1,\cdots, N-1,
\end{eqnarray}
where $l^{*}$, $g^{*}$ are dependent only on $l^{(1)}$, $l^{(2)}$, $l^{(3)}$, $g^{(1)}$, $g^{(2)}$ and the derivatives of $w^{(1)}_\tau$, $w^{(2)}_\tau$, $w_h$.

Then we define $\tilde{e}_i^n := W_i^n - x_{h_i}^n$, $0 \le i \le M$, $n=0,1,\cdots, N$. In other words, instead of a direct comparison between the numerical solution and exact PDE solution, we evaluate the numerical error between the numerical solution and the constructed solution $W$. The higher order truncation error enables us to obtain a required $W_h^{1,\infty}$ of the numerical solution, which is necessary in the nonlinear convergence analysis.

Note that the discrete $L^2$ norm $\|\tilde{e}^{0}\|_{2}=0$ at time step $t^0$. We make the following a-priori assumption at time step $t^n$:
\begin{equation}\label{priori}
\|\tilde{e}^{n}\|_2\leq  (\tau^{\frac{11}{4}}+ h^{\frac72} ) . 
\end{equation}
In turn,  the following estimates become available, by making use of inverse inequalities:
\begin{eqnarray}
&&
\|\widetilde{D}_h\tilde{e}^{n}\|_{2}\leq C (\tau^{\frac74} +h^{\frac52} ), \label{a-priori}\\
&&
 \|\widetilde{D}_h \tilde{e}^{n} \|_{\infty} \leq C C_m\frac{\|\widetilde{D}_h\tilde{e}^{n}\|_{2}}{h^{1/2}}\leq C C_m (\tau^{\frac54}+h^2), \ \mbox{if\ \ } h=O(\tau),
\label{a-priori_h} \\
&&
\|\widetilde{D}_h x^{n}_h \|_{\infty}=\|\widetilde{D}_h W^n-\widetilde{D}_h\tilde{e}^n\|_{\infty}\leq C^*+1 :=C_0^*, \label{a-priori_x}\\
&&
\mbox{with\ \ } C^*:=\|\widetilde{D}_h W^n\|_{\infty},\ \ \mbox{if\ \ } C C_m (\tau^{\frac54}+h^2 )\leq 1, \nonumber\\
&&
\Big\|\frac{\widetilde{D}_h x^{n}_h-\widetilde{D}_h x^{n-1}_h }{\tau}\Big\|_{\infty}=\Big\|\frac{\widetilde{D}_h W^{n}-\widetilde{D}_h W^{n-1} }{\tau}-\frac{\widetilde{D}_h \tilde{e}^{n}-\widetilde{D}_h \tilde{e}^{n-1} }{\tau}\Big\|_{\infty}\leq \tilde{C}^*_t+1,\label{a-priori_x_t}\\
&&
\mbox{with\ \ } \tilde{C}^*_t:=\Big\|\frac{\widetilde{D}_h W^{n}-\widetilde{D}_h W^{n-1} }{\tau}\Big\|_{\infty}, \ \mbox{if\ \ } C C_m (\tau^{\frac{1}{4}}+h )\leq 1. \label{a-priori-diff-time}
\end{eqnarray}
%
For $x_h,W\in\mathcal{Q}$, i.e., $\exists$ $\delta_0>0$, such that  $\widetilde{D}_h W^n_i\geq \delta_0$, then  $\widetilde{D}_h x^n_{h_i}\geq \frac{\delta_0}{2}>0$, $0\leq i\leq M$, if  $C_m \gamma(\tau^{\frac54}+h^2 )\leq\frac{\delta_0}{2}$.



In turn, subtracting \eqref{num:U} from the numerical scheme  \eqref{equ:numnum0} yields
\begin{eqnarray}
  && \ \ \ \frac{f_0(X_i)}{m\big(\frac{f_0(X)}{\widetilde{D}_h x^{n}_{h}}\big)_i^{m-1}}\cdot\frac{\tilde{e}^{n+1}_i - \tilde{e}^n_i}{\tau} +
  \frac{f_0(X_i)}{ m[f_0(X_i)]^{m-1}}\cdot\frac{W^{n+1}_i-W^{n}_i}{\tau}\cdot[\big(\widetilde{D}_h W^{n}\big)_i^{m-1}-\big(\widetilde{D}_h x_h^{n}\big)_i^{m-1}]
  \nonumber \\
  &&
  = d_h \left( \frac{f_0(X)}{D_h W_i^{n+1} D_h x_h^{n+1}}  D_h \tilde{e}^{n+1}\right)_i   + \tau^3 l^*_i+h^4g^{*}_i,\ \ 1\leq i\leq M-1,  \nonumber \\
  && \mbox{with\ }\tilde{e}^{n+1}_0 = \tilde{e}^{n+1}_M = 0,
  \label{consistency-3}
\end{eqnarray}
in which the form of the  left term comes from the following identity:
 \begin{eqnarray}
 &&
 \ \ \ \frac{f_{0}(X_i)}{m\big(\frac{f_0(X)}{\widetilde{D}_h W^{n}}\big)_i^{m-1}}\frac{W^{n+1}_i-W^{n}_i}{\tau}-
 \frac{f_{0}(X)}{m\big(\frac{f_0(X_i)}{\widetilde{D}_h x_h^{n}}\big)_i^{m-1}}\frac{x^{n+1}_{h_i}-x^{n}_{h_i}}{\tau}
  \nonumber\\
 &&
 = \frac{f_{0}(X_i)}{\tau m [f_0(X_i)]^{m-1}}[(\widetilde{D}_h W^{n})_i^{m-1}(W^{n+1}_i-W^{n}_i)
 -(\widetilde{D}_h x^{n}_{h})_i^{m-1}(x^{n+1}_{h_i}-x^{n}_{h_i}) \nonumber\\
 &&
 \  \  \ +\big(\widetilde{D}_h x^{n}_h\big)_i^{m-1}(W^{n+1}_i-W^{n}_i)-\big(\widetilde{D}_h x^{n}_h\big)_i^{m-1}(W^{n+1}_i-W^{n}_i)]\nonumber\\
 &&
 =\frac{f_0(X_i)}{ m[f_0(X_i)]^{m-1}}\cdot\frac{W^{n+1}_i-W^{n}_i}{\tau}\cdot[\big(\widetilde{D}_h W^{n}\big)_i^{m-1}-\big(\widetilde{D}_h x_h^{n}\big)_i^{m-1}] \nonumber \\
 &&
 \ \ \ +\frac{f_0(X_i)}{m\big(\frac{f_0(X)}{\widetilde{D}_h x_h^{n}}\big)_i^{m-1}}  \cdot\frac{\tilde{e}_i^{n+1} - \tilde{e}_i^n}{\tau}. \nonumber
  \end{eqnarray}

Based on the preliminary results, taking a discrete inner product with (\ref{consistency-3}) by $2 \tilde{e}^{n+1}$ gives
\begin{eqnarray}
&&
 2\left\langle \alpha_n (\tilde{e}^{n+1} - \tilde{e}^n), \tilde{e}^{n+1}\right\rangle -2\tau\left\langle  d_h \left( \frac{f_0(X)}{D_h W^{n+1} D_h x^{n+1}_h}  D_h \tilde{e}^{n+1}\right), \tilde{e}^{n+1}\right\rangle\nonumber\\
 &&
 =-2\tau\left\langle \frac{f_0(X)}{ m[f_0(X)]^{m-1}}\cdot\frac{W^{n+1}-W^{n}}{\tau}\cdot[(\widetilde{D}_h W^{n})^{m-1}-(\widetilde{D}_h x^{n}_h)^{m-1}], \tilde{e}^{n+1}\right\rangle \nonumber\\
 &&
 \ \ \ +2\tau\left\langle \tau^3 f^*+h^4g^{*},\tilde{e}^{n+1}\right\rangle,
  \label{convergence-1}
\end{eqnarray}
where  \begin{equation}\label{alpha}\alpha_n:=\frac{f_0(x)}{m\big(\frac{f_0(X)}{\widetilde{D}_h x^{n}_h}\big)^{m-1}}.\end{equation}

For the first  term of the left side, we get
\begin{equation}
\begin{split}
2\left\langle \alpha_n  (\tilde{e}^{n+1} - \tilde{e}^n),\tilde{e}^{n+1}\right\rangle
&
=\alpha_n \|\tilde{e}^{n+1}\|^2_{2}+\alpha_n \|\tilde{e}^{n+1}-\tilde{e}^{n}\|^2_{2}-\alpha_n \|\tilde{e}^{n}\|^2_{2} \\
&
\geq \alpha_n \|\tilde{e}^{n+1}\|^2_{2}-\alpha_n \|\tilde{e}^{n}\|^2_{2}. \label{convergence-proof-1}
\end{split}
\end{equation}

For the second term of the left side, we see that
\begin{eqnarray}
&&
 -2\tau\left\langle  d_h \left( \frac{f_0(X)}{D_h W^{n+1} D_h x^{n+1}_h}  D_h \tilde{e}^{n+1}\right), \tilde{e}^{n+1}\right\rangle \nonumber\\
&&
=2\tau\left\langle \frac{f_0(X)}{D_h W^{n+1} D_h x_h^{n+1} }  D_h \tilde{e}^{n+1}, D_h \tilde{e}^{n+1}\right\rangle_e \geq 0,
%
\label{convergence-proof-2} 
\end{eqnarray}
in which  the summation by parts formula (\ref{FD-inner product-3-1}) is applied with $\tilde{e}^{n+1}_0 = \tilde{e}^{n+1}_N = 0$.

For the right side term, we have
\begin{eqnarray}
&&
-2\tau\left\langle \frac{f_0(X)}{ m[f_0(X)]^{m-1}}\cdot\frac{W^{n+1}-W^{n}}{\tau}\cdot[(\widetilde{D}_h W^{n})^{m-1}-(\widetilde{D}_h x^{n}_h)^{m-1}], \tilde{e}^{n+1}\right\rangle \nonumber\\
&&
=-2\tau\left\langle \frac{f_0(X)}{m[f_0(X)]^{m-1}}\cdot\frac{W^{n+1}-W^{n}}{\tau}\cdot[(m-1)(\widetilde{D}_h \zeta^{n})^{m-2}\widetilde{D}_h \tilde{e}^{n}], \tilde{e}^{n+1}\right\rangle \nonumber\\
&&
\leq 2\tau C_1\|\widetilde{D}_h  \tilde{e}^{n}\|_2\|\tilde{e}^{n+1}\|_2, \ \ \Big(C_1:=\frac{(m-1)B_fC^*_tC_{\zeta}^{m-2}}{mb_f^{(m-1)}} \Big),\\
&&
\leq\tau C_1\|\widetilde{D}_h  \tilde{e}^{n}\|_2^2+\tau C_1\|\tilde{e}^{n+1}\|_2^2, \nonumber
 \label{convergence-proof-3}
\end{eqnarray}
in which $ C^*_t=\|W_t\|_{\infty}$,
$\widetilde{D}_h\zeta$ is between $\widetilde{D}_h x_h^n$, and $\widetilde{D}_h W^n$,  $\|\widetilde{D}_h\zeta\|_{\infty}\leq C_\zeta$, with
\begin{equation*}
C_\zeta:=
\begin{cases}
C^*_0, & m\geq2,\nonumber\\
\delta_0/2, & m<2.
\end{cases}
\end{equation*}

The local truncation error term could be bounded by the standard Caught inequality:
\begin{equation}
\begin{split}
2\tau\left\langle \tau^3 l^*+h^4g^{*}, \tilde{e}^{n+1}\right\rangle
&
\leq \tau\|\tau^3 l^*+h^4g^{*}\|_2^2+\tau\|\tilde{e}^{n+1}\|^2_2  \\
&
\leq \tau C(\tau^3+h^4)^2+\tau\|\tilde{e}^{n+1}\|^2_2.
  \label{convergence-proof-4}
  \end{split}
\end{equation}



Next we estimate $\|D_h x^{n+1}_h\|_{\infty}$ roughly.
Based on \eqref{alpha}, $\alpha_n$ can be estimated by
$$C_\alpha:=\frac{b_f}{m(\frac{B_f}{\delta_0/2})^{m-1}}\leq\|\alpha\|_n\leq\frac{B_f}{mb_f^{m-1}}(C_0^*)^{m-1}:=\bar{C}_\alpha.$$
A substitution of \eqref{convergence-proof-1}-\eqref{convergence-proof-4} into (\ref{convergence-1}), in combination with \eqref{a-priori}, leads to
\begin{equation*}
\begin{split}
 (\alpha_n-\tau (1+C_1))\| \tilde{e}^{n+1} \|_{2}^2
 &
 \leq \alpha_n\| \tilde{e}^n \|_{2}^2+\tau C_1 \|\widetilde{D}_h  \tilde{e}^{n}\|_2^2+\tau C(\tau^3+h^4)^2 \nonumber\\
 &
 \leq \tau \bar{C}(\tau^\frac74 + h^\frac52)^2,
  \end{split}
\end{equation*}
where $\bar{C}$ is dependent on $C$, $C_1$ and $\bar{C}_\alpha$.
Then we get
\begin{eqnarray}
&&
  \| \tilde{e}^{n+1} \|_{2}^2 \leq \widetilde{C}^2\tau(\tau^\frac74+h^\frac52)^2, 
  \mbox{\ \ i.e.,\ \ } \| \tilde{e}^{n+1} \|_{2} \leq \widetilde{C}\tau^{\frac{1}{2}}(\tau^\frac74 + h^\frac52), \end{eqnarray}
  \mbox{with\ \ } $\widetilde{C}:=\Big(\frac{\bar{C}}{C_\alpha/2}\Big)^{\frac{1}{2}}$, \ \ \
 \mbox{if\ \ } $\tau (1+C_1)\leq C_\alpha/2$.

Based on  the inverse inequality \eqref{equ:inverse}, we obtain that, by choosing  $h=O(\tau)$,
\begin{equation}
\| \tilde{e}^{n+1} \|_{\infty}\leq \frac{C_m\| \tilde{e}^{n+1} \|_{2}}{h^{\frac{1}{2}}}\leq C_m\widetilde{C}(\tau^\frac74 + h^\frac52).
\end{equation}
Then we have
\begin{eqnarray}
\|D_h {x}^{n+1}_h \|_{\infty}=\|D_h {W}^{n+1} -D_h \tilde{e}^{n+1} \|_{\infty}\leq C^*+C_m\widetilde{C}(\tau^\frac34 +h^\frac32)\leq C^*+1:=C^*_0, \label{rough_x}\ \ \
\end{eqnarray}
if $C_m\widetilde{C}(\tau+h^2)\leq 1$.

As a result, \eqref{convergence-proof-2} can be re-estimated as follows:
\begin{equation}\label{convergence-proof-2-2}
2\tau\left\langle \frac{f_0(X)}{D_h W^{n+1} D_h x_h^{n+1}}  D_h \tilde{e}^{n+1}, D_h \tilde{e}^{n+1}\right\rangle_e
\geq 2\tau C_2\|D_h \tilde{e}^{n+1}\|_2^2,
\end{equation}
with $C_2:=\frac{b_f}{C^*C_0^*}$.

As a consequence, a substitution of \eqref{convergence-proof-1}-\eqref{convergence-proof-4} with \eqref{convergence-proof-2-2} into (\ref{convergence-1}) leads to
\begin{eqnarray}
\alpha_n\| \tilde{e}^{n+1} \|_{2}^2 - \alpha_n\| \tilde{e}^n \|_{2}^2   +
 \tau C_2\|\widetilde{D}_h \tilde{e}^{n+1}\|^2_{2} 
 \leq \tau\Big(1+\frac{C_1^2}{C_2}\Big)\|\tilde{e}^{n+1}\|^2_{2}
 +\tau C (\tau^3+h^4)^2,\nonumber
\end{eqnarray}
where the following estimates are applied: $\|\widetilde{D}_h x^n_h\|_2\leq\|D_h x^n_h\|_2$ and
\begin{equation}
2\tau C_1\|\widetilde{D}_h  \tilde{e}^{n}\|_2\|\tilde{e}^{n+1}\|_2\leq
\tau\frac{C_1^2}{C_2}||\tilde{e}^{n+1}||^2_{2} +\tau C_2||\widetilde{D}_h \tilde{e}^{n}||^2_{2} .
\end{equation}

Subsequently, a summation in time shows that 
\begin{eqnarray}
 \alpha_n \| \tilde{e}^{n+1} \|_{2}^2+
  \tau C_2\sum\limits_{k=1}^{n+1}||\widetilde{D}_h \tilde{e}^{k}||^2_{2}
 & \leq &\tau\sum\limits_{k=1}^{n}\frac{(\alpha_k- \alpha_{k-1})}{\tau}\|\tilde{e}^{k}\|^2_2+
  \tau \Big(\frac{C_1^2}{C_2}+1\Big)\sum\limits_{k=1}^{n+1}\|\tilde{e}^{k}\|^2_{2} \notag\\
  && + CT (\tau^3+h^4)^2 , \notag \\
\| \tilde{e}^{n+1} \|_{2}^2+
  \tau \frac{C_2}{C_\alpha}\sum\limits_{k=1}^{n+1}||\widetilde{D}_h \tilde{e}^{k}||^2_{2}
& \leq&  \frac{\tau}{C_\alpha} (\frac{C_1^2}{C_2}+1+\widetilde{C}_{\alpha})\sum\limits_{k=1}^{n+1}\|\tilde{e}^{k}\|^2_{2} +\frac{CT}{C_\alpha} (\tau^3+h^4)^2,\notag
\end{eqnarray}
where we have used the estimate
\begin{equation}
\begin{split}
\big\|\frac{\alpha^k-\alpha^{k-1}}{\tau}\big\|_{\infty}
&
=\Big\|\frac{f_0(X)}{m[f_0(X)]^{m-1}}\cdot\frac{(\widetilde{D}_h x^{k}_h)^{m-1}-(\widetilde{D}_h x^{k-1}_h)^{m-1}}{\tau}\Big\|_{\infty} \nonumber\\
&
=\Big\|\frac{f_0(X)}{m[f_0(X)]^{m-1}}(m-1)(\widetilde{D}_h \vartheta)^{m-2}\frac{\widetilde{D}_h x^{k}_h-\widetilde{D}_h x^{k-1}_h}{\tau}\Big\|_{\infty}\nonumber\\
&
\leq \frac{(m-1)B_f}{m b_f^{m-1}}(C_{\vartheta})^{m-2}(\tilde{C}^*_t+1):=\widetilde{C}_{\alpha} . \nonumber
 \end{split}
\end{equation}
It is noticed that $T$ is the terminal time, \eqref{a-priori-diff-time} is applied and  $\widetilde{D}_h \vartheta$ is between $\widetilde{D}_h x^k_h$ and $\widetilde{D}_h x^{k-1}_h$ with
\begin{equation*}
\|\widetilde{D}_h \vartheta\|_{\infty}\leq C_{\vartheta}:=
\begin{cases}
C^*_0, & m\geq2,\\
\delta_0/2, & m<2.
\end{cases}
\end{equation*}

In turn, an application of discrete Gronwall inequality yields the desired convergence result:
\begin{eqnarray}
  \| \tilde{e}^{n+1} \|_{2}^2   +
  \tau \frac{C_2}{C_\alpha}\sum\limits_{k=1}^{n+1}||\widetilde{D}_h \tilde{e}^{k}||^2_{2}
  \leq e^{TC_0}\frac{CT}{C_\alpha}(\tau^3+h^4)^2, \quad
 \mbox{i.e.,} \, \,  \| \tilde{e}^{n+1} \|_{2}\leq \gamma(\tau^3+h^4), \nonumber
\end{eqnarray}
where  $C_0:=\frac{1}{C_\alpha} (\frac{C_1^2}{C_2}+\widetilde{C}_{\alpha}+1)$ and
\begin{equation}
\gamma:=\Big(\frac{CT }{C_\alpha}\Big)^{\frac{1}{2}}e^{\frac{C_0T}{2}}.\label{M}
\end{equation}
Therefore, the a-priori assumption \eqref{priori} is also valid at $t^{n+1}$:
\begin{equation}
\| \tilde{e}^{n+1} \|_{2}\leq \gamma(\tau^3+h^4) \le \tau^{\frac{11}{4}} + h^{\frac72} ,
\end{equation}
provided that $\tau \leq \gamma^{-4}$, $h \le \gamma^{-2}$. 

 Based on the following estimate
  \begin{equation}\label{derivative}
 \|\widetilde{D}_h \tilde{e}^{n+1}\|_2 =\|\widetilde{D}_h x^{n+1}_h-\widetilde{D}_h W^{n+1}\|_2\leq C \gamma(\tau^2+h^3),
 \end{equation}
 we obtain
  \begin{equation}\label{derivative_u_e}
  \|\widetilde{D}_h x_h^{n+1}-  \widetilde{D}_h x_e^{n+1}\|_2\leq C(\tau+h^2).
  \end{equation}
Finally, we estimate the error  between the numerical solution $f_h^{n+1}$ and the exact solution $f_e^{n+1}$ of the problem \eqref{eqcm}-\eqref{eqbc}:
 \begin{equation}
 \begin{split}
\|f_{e}^{n+1}-f_{h}^{n+1}\|_2&=\Big\|\frac{f_0(X)}{\partial_X x^{n+1}_e}-\frac{f_0(X)}{\widetilde{D}_h x^{n+1}_h}\Big\|_2 \nonumber \\
&
=\left\|\frac{f_0(X)}{\partial_X x^{n+1}_e}-\frac{f_0(X)}{\widetilde{D}_h x^{n+1}_e}+\frac{f_0(X)}{\widetilde{D}_h x^{n+1}_e}-\frac{f_0(X)}{\widetilde{D}_h x^{n+1}_h}\right\|_2  \nonumber \\
&
\leq C(\tau+h^2).
\end{split}
\end{equation}
 $\hfill\Box$
\section{Numerical Results}
\label{sec:4}
In this section, we present some numerical results to demonstrate the convergence rate of the numerical scheme.

Before that, we define the error of a numerical solution measured  in the $\mathcal{L}^2$ and $\mathcal{L}^{\infty}$ norms as:
\begin{equation}\label{L2}
\|e_h\|_2^2=\frac{1}{2}\left(e_{h_0}^2 h_{x_{0}}+\sum\limits_{i=1}^{M-1}e_{h_i}^2 h_{x_i}+e_{h_M}^2 h_{x_M}\right),
\end{equation}
and
\begin{equation}\label{Linf}
\|e_h\|_{\infty}=\max\limits_{0\leq i\leq M}\{|e_{h_i}|\},
\end{equation}
 where $e_h=(e_{h_0},e_{h_1},\cdots,e_{h_M})$ and
for the error of the density  $f-f_h$,
 \begin{equation*}
 h_{x_i}=x_{i+1}-x_{i-1}, \ \ 1\leq i \leq M-1; \ \ \
 h_{x_0}=x_{1}-x_{0}; \ \
 h_{x_M}=x_{M}-x_{M-1}, \end{equation*}
and for  the error of the  trajectory $x-x_h$,
 \begin{equation}
h_{x_i}=2h,  \ \ 1\leq i\leq M-1, \ \ \ h_{x_0}=h_{x_M}=h , \nonumber
\end{equation}
 where $h$ is the spatial step.

 Consider the problem \eqref{eqcm}-\eqref{eqbc} in dimension one with a smooth positive initial data
\begin{equation}\label{eqiniEx}
f_0(x)=\frac{1}{2}(-x^2+1.01), x\in\Omega:=[-1,1].
\end{equation}

Firstly, the trajectory equation \eqref{eqtra} with the initial and boundary condition \eqref{eqtrabou}-\eqref{eqtraini} can be solved by the fully discrete scheme \eqref{equ:numnum0}. And then the density function $f$ in \eqref{equ:conservationL} can be approximated by \eqref{eqDen}. The reference ``exact" solution is obtained numerically on a much finer mesh with $h=\frac{1}{10000},\  \tau=\frac{1}{10000}$.

\begin{figure}
\captionsetup{font={scriptsize}}
\centering
\subfigure[\scriptsize $t=0.1$]
{\includegraphics[width=6cm,height=5cm]{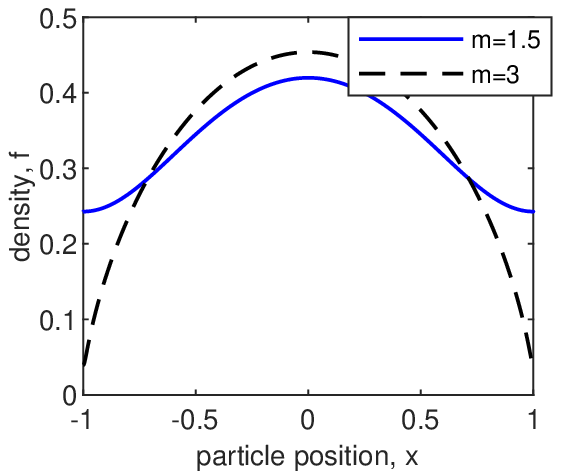}}\hspace{.1in}
\subfigure[\scriptsize $t=0.5$]
{\includegraphics[width=6cm,height=5cm]{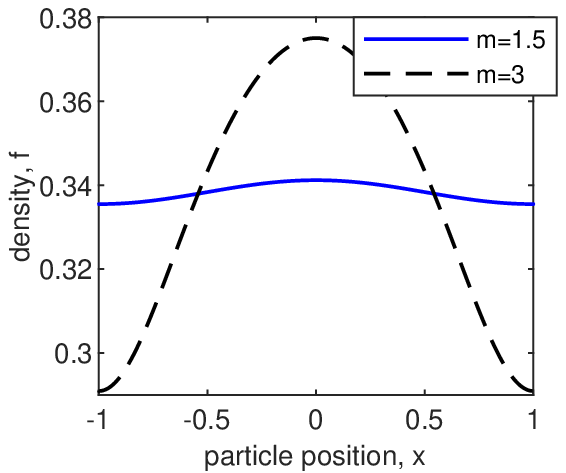}}\hspace{.1in}
\caption{The evolution of $f$ ($h=1/1000$, $\tau=1/1000$)}
\label{fig:Den}
\end{figure}

\begin{figure}
\begin{minipage}[t]{0.45\linewidth}
\captionsetup{font={scriptsize}}
\centering
\includegraphics[width=6cm,height=5cm]{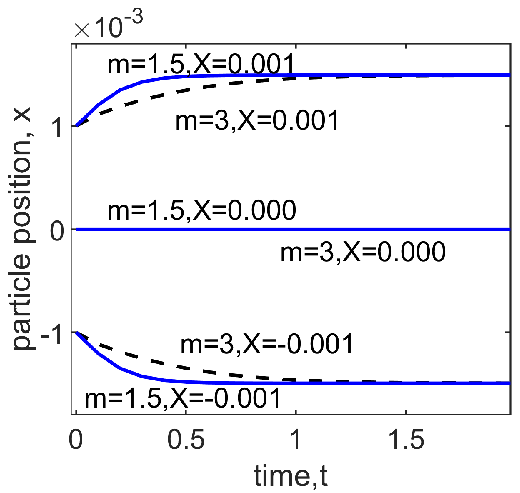}
\caption{The evolution of particle with initial \\ position $X= -0.001, 0.000, 0.001$ ($h=1/1000$, \\$\tau=1/1000$)}
\label{fig:Particle}
\end{minipage}
\begin{minipage}[t]{0.45\linewidth}
\captionsetup{font={scriptsize}}
\centering
\includegraphics[width=6cm,height=5cm]{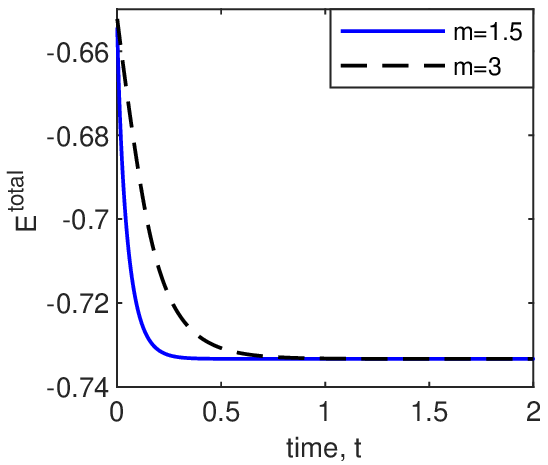}
\caption{The evolution of total energy ($h=1/1000$, $\tau=1/1000$)}
\label{fig:Energy}
\end{minipage}
\end{figure}

Table \ref{tablePME11}  shows  the convergence rate with $m=1.5$ and $m=3$ at time $T=0.5$. The rate for  density $f$ and trajectory $x$ in the   $\mathcal{L}^2$ and  $\mathcal{L}^{\infty}$ norm is  second order in space and first order in time without dependence on $m$.  Fig. \ref{fig:Den} presents the  density $f$ at time $t=0.1$ and $t=0.5$ for both values of $m$. The results imply that the speed of diffusion  decreases as $m$ increases. Fig. \ref{fig:Particle} displays the evolution of particles whose initial positions are $X=-0.001,\ 0.000,\ 0.001$, respectively. We see that particles move outward at a finite speed. However, the speed is lower as $m$ increases, except for the center point which remains stationary.  As shown in Fig. \ref{fig:Energy},  the total energy decays as time evolves for both values of $m$, and the decreasing rate is slowed down as $m$ increases.

\begin{threeparttable}[b]
\scriptsize
\centering
\setlength{\abovecaptionskip}{10pt}
\setlength{\belowcaptionskip}{-3pt}
\caption{\scriptsize  Convergence  rate of solution $f$ and trajectory $x$ at time $T=0.5$}
\label{tablePME11}
\begin{tabular}{@{ } l c c c c c c c c c }

\hline
\multicolumn{1}{l}{}
&\multicolumn{8}{c}{$m=1.5$}\\\hline
 $h$    &$\tau$ &$ \mathcal{L}^2$-error $ (f) $    & Order   &$ \mathcal{L}^{\infty}$-error$ (f) $    & Order &$ \mathcal{L}^2$-error $(x)$   & Order  &$ \mathcal{L}^{\infty}$-error $(x)$   & Order  \\\hline
1/10 &1/10 & 1.5683e-02&  &1.1680e-01& &1.5122e-03 & &5.2129e-03& \\\hline
 1/20 &1/40 &2.8389e-03 & 2.7621&3.8639e-02&1.5114&1.4812e-03 &2.3315&1.1415e-03&2.2833\\\hline
 1/40 &1/160 &4.8401e-04 &2.9327 &1.0515e-02&1.8373&3.1765e-04 &2.1501&2.6587e-04&2.1468\\\hline
 1/80 &1/640 &9.2853e-05 &2.6063 &2.6850e-03&1.9582&7.3867e-05 &2.1474&6.1999e-05&2.1441\\\hline

\hline
\multicolumn{1}{l}{}
&\multicolumn{8}{c}{$m=3$}\\\hline
  $h$    &$\tau$ &$ \mathcal{L}^2$-error $ (f) $    & Order   &$ \mathcal{L}^{\infty}$-error$ (f) $    & Order &$ \mathcal{L}^2$-error $(x)$   & Order  &$ \mathcal{L}^{\infty}$-error $(x)$   & Order \\\hline
1/10 &1/10 & 1.7044e-02&  &1.1471e-01& &1.9900e-03 & &7.1606e-03& \\\hline
1/20 &1/40 &3.0278e-03 & 2.8146&3.5746e-02&1.6045&4.9937e-04 &1.9924&1.8319e-03&1.6045\\\hline
1/40 &1/160 &5.5190e-04 &2.7431 &9.6637e-03&1.8495&1.2332e-04 &2.0247&4.5340e-04&1.8495\\\hline
1/80 &1/640 &1.1139e-04 &2.4773 &2.4584e-03&1.9655&2.9237e-05 &2.1089&1.0757e-04&1.9655\\\hline
\end{tabular}
\begin{tablenotes}
     \scriptsize
        \item[1]   $ \mathcal{L}^2$-error and $ \mathcal{L}^{\infty}$-error  is defined by  \eqref{L2} and \eqref{Linf}, respectively.

                \item[2]    $\tau$ is the time step and  $h$ is the space step.
        \end{tablenotes}
\end{threeparttable}

More interesting examples can be found in \cite{C. Duan(2018)}, such as a free boundary problem with a exact Barenblatt solution, the waiting time phenomenon and the problem with two support sets at the initial state.
 \section{Conclusion}
  The numerical methods of the PME based on EnVarA has been proposed and studied in \cite{C. Duan(2018)}, while a theoretical justification for optimal convergence analysis has not been available. In this paper, we prove the second order spatial convergence and the first order temporal convergence for the nonlinear numerical scheme.  A careful asymptotic expansion for the exact solution in terms of the numerical scheme is applied to obtain higher order consistency. Furthermore, we use two step error estimates: a rough estimate to control a discrete $W^{1,\infty}$ bound of the numerical solution, and a refined estimate to derive the desired convergence result. 

One obvious limitation of this work is associated with the one-dimensional nature of the problem. In two or higher dimension, the determinant of the deformation gradient, i.e., $\det\frac{\partial x}{\partial X}$, will arise in the trajectory equation, which is a  complex nonlinear degenerate parabolic equation system. A suitable numerical method in multi-dimensional case, which can satisfy the discrete energy dissipation law, is still in the investigation process. Solving for multi-dimensional PME by this energetic method and the corresponding optimal error estimate will be left to the future works. Another limitation is that the assumption of a positive initial condition ($f_0>0$), in which the convergence rate does not depend on the constant $m$. It is well known that if the initial state has a compact support, the convergent rate decreases with  $m$. In this case, the trajectory equation with a free boundary makes the convergence analysis more difficult. This problem will also be considered in the future works.

\noindent\emph{Acknowledgments.}
The work of Yue is supported in part by NSF of China under the grants 11971342. Chun Liu and Cheng Wang are partially supported by NSF grants DMS-1216938, DMS-1418689, respectively.

 {\footnotesize
\bibliographystyle{unsrt}

}

 \end{document}